\documentclass[a4paper]{article}
\usepackage{amssymb,amsmath,amsthm}
\usepackage{graphicx}
\usepackage{hyperref}
\newcommand{\ie}{\emph{i.e.}}
\newcommand{\eg}{\emph{e.g.}}
\newcommand{\cf}{\emph{cf}}
\newcommand{\Nat}{\mathbb{N}}
\newcommand{\Real}{\mathbb{R}}

\newcommand{\Int}{\mathbb{Z}}
\newcommand{\sii}{L^2}
\newcommand{\Dom}{\mathfrak{D}}

\newcommand{\rot}{\mathop{\mathrm{rot}}\nolimits}
\newcommand{\dist}{\mathop{\mathrm{dist}}\nolimits}
\newcommand{\eps}{\varepsilon}
\newcommand{\const}{\mathrm{const}}

\newtheorem{Theorem}{Theorem}
\newtheorem{Proposition}{Proposition}
\newtheorem{Corollary}{Corollary}
\theoremstyle{remark}
\newtheorem{Remark}{Remark}
\begin{document}
%
\title{\textbf{\Large
The improved decay rate for the heat semigroup
with local magnetic field in the plane
}}
\author{
David Krej\v{c}i\v{r}{\'\i}k\,%
\footnote{On leave from
\emph{Department of Theoretical Physics,
Nuclear Physics Institute ASCR,
25068 \v Re\v z, Czech Republic; 
krejcirik@ujf.cas.cz}.}
}
\date{\small
%
\emph{Basque Center for Applied Mathematics, 
Bizkaia Technology Park, \\ Building 500, 
48160 Derio, Kingdom of Spain}
\medskip \\
\emph{IKERBASQUE, Basque Foundation for Science, \\ 
48011 Bilbao, Kingdom of Spain}
%
\ \bigskip \\
10 January 2011
}
\maketitle
\begin{abstract}
\noindent
We consider the heat equation in the presence
of compactly supported magnetic field in the plane. 
We show that the magnetic field
leads to an improvement of the decay rate of 
the heat semigroup by a polynomial factor 
with power proportional to the distance 
of the total magnetic flux to the discrete set of flux quanta. 

The proof employs Hardy-type inequalities due to Laptev and Weidl
for the two-dimensional magnetic Schr\"odinger operator 
and the method of self-similar variables
and weighted Sobolev spaces for the heat equation.
A careful analysis of the asymptotic behaviour 
of the heat equation in the similarity variables
shows that the magnetic field asymptotically degenerates to 
an Aharonov-Bohm magnetic field with the same total magnetic flux,
which leads asymptotically to the gain on the polynomial decay rate
in the original physical variables.

Since no assumptions about the symmetry of the magnetic field
are made in the present work, it confirms that 
the recent results of Kova\v{r}{\'\i}k~\cite{Kovarik_2010}
about large-time asymptotics of the heat kernel of 
magnetic Schr\"o\-din\-ger operators
with radially symmetric field hold in greater generality. 

%
%
\end{abstract}
%
%
\newpage
\section{Introduction}
%
The principal objective of this paper is to show that
the solutions of the heat equation
\begin{equation}\label{heat.intro}
  u_t + \big(\!-i\nabla_{\!x}-A(x)\big)^2 u = 0
  \,,
\end{equation}
in space-time variables $(x,t) \in \Real^2\times(0,\infty)$,
admit a faster large-time decay due to
the presence of a magnetic field $B = \rot A$
that is supposed to vanish at infinity.
The characteristic property under which we prove the result
is that the \emph{total magnetic flux}
\begin{equation}\label{flux}
  \Phi_B := \frac{1}{2\pi} \int_{\Real^2} B(x) \, dx
\end{equation}
does not belong to the discrete set of \emph{flux quanta},
which coincides with integers for our choice of physical constants.
As a measure of the additional decay of solutions,
we consider the (polynomial) \emph{decay rate}
\begin{equation}\label{rate}
  \gamma_B
  := \sup \Big\{ \gamma \left| \
  \exists C_\gamma > 0, \, \forall t \geq 0, \
  \big\|e^{-(-i\nabla-A)^2 t}\big\|_{
  \sii(K)
  \to
  \sii
  }
  \leq C_\gamma \, (1+t)^{-\gamma}
  \Big\} \right.
  .
\end{equation}
Here we use the shorthands
$\sii:=\sii(\Real^2)$ and $\sii(K):=\sii(\Real^2,K(x)\,dx)$
with the Gaussian weight
\begin{equation}\label{weight}
  K(x) := e^{|x|^2/4}
  \,.
\end{equation}
Our main result reads as follows:
\begin{Theorem}\label{Thm.rate}
Let $B \in C_0^0(\Real^2)$.
We have
$$
  \gamma_B
  \begin{cases}
    \, = 1/2
    & \mbox{if $B=0$},
    \\
    \, \geq (1+\beta)/2
    & \mbox{otherwise},
  \end{cases}
$$
where $\beta:=\dist(\Phi_B,\Int)$.
\end{Theorem}

Here the result for $B=0$ follows easily by the explicit knowledge
of the heat kernel for the Laplacian in~$\Real^2$.
The essential content of Theorem~\ref{Thm.rate}
is that solutions to~\eqref{heat.intro}
when the total flux~$\Phi_B$ is not an integer
gain a further polynomial decay rate $\beta/2$.
The proof of this statement is more involved
and constitutes the main body of the paper.
It is based on the method of self-similar solutions
recently developed for twisted-waveguide systems
by Zuazua and the present author in~\cite{KZ1,KZ2}.

The result has a rather natural interpretation in terms of
\emph{diamagnetism} that we briefly explain now.

\subsection{Relation to Hardy inequalities}
Diamagnetism in quantum mechanics is reflected in
repulsive properties of the Hamiltonian
\begin{equation}\label{Hamiltonian}
  H_B := (-i\nabla-A)^2
\end{equation}
describing a spin-less particle in the external magnetic field~$B$.
This is most easily seen for the homogeneous field $B(x)=B_0$
for which the ground-state energy is raised by the field magnitude,
\ie\ $\inf\sigma(H_{B_0})=|B_0|$,
and any such non-zero field therefore leads to
a Poincar\'e-type inequality $H_{B_0} \geq |B_0|$.
The latter remains true for magnetic fields
satisfying $\pm B(x) \geq |B_0|$.
This clearly implies an exponential decay rate
$
  \|e^{- H_B t}\|_{\sii\to\sii} \leq e^{-|B_0| t}
$
in view of the spectral mapping theorem.

For magnetic fields vanishing at infinity considered in this paper,
the repulsive effect is more subtle because the spectrum starts by zero:
$
  \sigma(H_{B}) = [0,\infty)
$,
regardless of the presence of the magnetic field.
Hence,
$
  \|e^{- H_B t}\|_{\sii\to\sii} = 1
$.
However, Laptev and Weidl proved in~\cite{Laptev-Weidl_1999}
an important Hardy-type inequality
\begin{equation}\label{Hardy}
  H_B \geq \frac{c}{1+|\cdot|^2}
\end{equation}
in the sense of quadratic forms, 
where $c=c(A)$ is a positive constant whenever $\Phi_B \not\in \Int$
(for generalizations, see 
\cite{Weidl_1999, Balinsky-Laptev-Sobolev_2004, Evans-Lewis_2005}).
On the other hand, because of the criticality of the free Hamiltonian in~$\Real^2$,
\eqref{Hardy}~cannot hold for $A=0$.
It is natural to expect
-- and supported by the general conjectures raised in~\cite{KZ1,FKP} --
that the additional repulsiveness caused by the Hardy-type inequality
will lead to an improvement of the asymptotic decay rate
of solutions to~\eqref{heat.intro},
provided that the initial data are restricted to a subspace of~$\sii$.
In fact, \eqref{Hardy}~plays a central role
in our proof of Theorem~\ref{Thm.rate}.

For dimensions greater or equal to three,
the inequality \eqref{Hardy}~is just a direct consequence
of the diamagnetic inequality~\eqref{diamagnetic}
and the classical Hardy inequality
for the free Hamiltonian,
so that~\eqref{Hardy} holds with a positive constant independent of~$A$
in the higher dimensions.
That is why the most interesting situation
is the two-dimensional case,
to which we restrict in this paper.

\subsection{Relation to previous results of Kova\v{r}\'ik \cite{Kovarik_2010}}
Before elaborating on the details of the proof of Theorem~\ref{Thm.rate},
let us make a few comments on the novelty of the results.

The present paper is partially motivated by a recent study
of Kova\v{r}\'ik \cite{Kovarik_2010} about large-time
properties of the heat kernel of two-dimensional
Schr\"odinger (and Pauli) operators
with \emph{radially symmetric} magnetic field.
Among other things, the author of~\cite{Kovarik_2010} shows
that the magnetic heat kernel has a faster decay for large times
than the heat kernel of the free Hamiltonian.
In some cases, he even shows that the additional decay rate
is exactly that of Theorem~\ref{Thm.rate}.
His results are more precise in the sense that
he gets pointwise (sometimes two-sided) estimates
and cover some cases with $\Phi_B \in \Int$.

On the other hand, the method of~\cite{Kovarik_2010}
is based on partial-wave decomposition and
is therefore  restricted to radially symmetric fields.
The method of the present paper, on the contrary, is more general
and it is actually relying only on the existence of
the magnetic Hardy-type inequality~\eqref{Hardy}.
In this way, we have been able to treat the general magnetic field
without any symmetry restriction.

Finally, Kova\v{r}\'ik calculated in~\cite{Kovarik_2010}
the exact large-time asymptotics for the Aha\-ro\-nov-Bohm field.
These asymptotics yield the same decay rates
as for regular magnetic fields, suggesting a somewhat universal
nature of the decay rates.
The present paper provides an insight into the phenomenon:
\emph{all the magnetic fields effectively behave for large times
as an Aharonov-Bohm field}.
This statement will be made precise in the following subsection.

\subsection{The idea of the proof}
As already mentioned above,
the principal idea behind the main result of Theorem~\ref{Thm.rate},
\ie\ the better decay rate due to the presence of magnetic field,
is the positivity of the right hand side of~\eqref{Hardy}.
In fact, Hardy-type inequalities have already been used
as an essential tool to study the asymptotic behaviour of the heat equation
in other situations
\cite{Cabre-Martel_1999,Vazquez-Zuazua_2000,KZ1,KZ2}.
However, it should be stressed that Theorem~\ref{Thm.rate} does not follow
as a direct consequence of~\eqref{Hardy} by some energy estimates
but that important and further
technical developments that we explain now are needed.
Nevertheless, overall, the main result of the paper again confirms
that the Hardy-type inequalities end up enhancing the decay rate
of solutions.

Let us now informally describe our proof
that there is the extra decay rate if the magnetic field is present.
\smallskip \\
\textbf{I (changing the space-time).}
The main ingredient in the proof is the method of self-similar solutions
developed in the whole Euclidean space
by Escobedo and Kavian~\cite{Escobedo-Kavian_1987}.
We rather follow the approach of \cite{KZ1,KZ2}
where the technique was adapted to twisted-waveguide systems
which exhibit certain similarities with the present problem.
Writing
\begin{equation}\label{SST}
  \tilde{u}(y,s) = e^{s/2} u(e^{s/2}y,e^s-1)
  \,,
\end{equation}
the equation~\eqref{heat.intro} is transformed to
\begin{equation}\label{heat.similar}
  \tilde{u}_s
  - \mbox{$\frac{1}{2}$} \, y \cdot \nabla_{\!y} \tilde{u}
  + \big(\!-i\nabla_{\!y}-A_s(y)\big)^2 \tilde{u}
  - \mbox{$\frac{1}{2}$} \, \tilde{u}
  = 0
\end{equation}
in ``self-similarity'' variables $(y,s)\in\Real^2\times(0,\infty)$,
where
\begin{equation}\label{sigma}
  A_s(y) := e^{s/2} A(e^{s/2}y)
  \,.
\end{equation}
When evolution is posed in that context,
$y$~plays the role of new space variable and~$s$ is the new time.
However, note that now we deal with a non-autonomous system.
\smallskip \\
\textbf{II (choosing rapidly decaying initial data).}
We reconsider~\eqref{heat.similar} in the weight\-ed space $\sii(K)$.
It has the advantage that the associated ($s$-dependent)
generator has compact resolvent in this space
and therefore purely discrete spectrum.
The latter is particularly useful for energy estimates.
\smallskip \\
\textbf{III (removing the weight).}
Just for personal reasons, we prefer to work
in the usual $\sii$-space, without the weight~$K$.
To achieve this, we make the substitution $v=K^{1/2}\tilde{u}$
and transform~\eqref{heat.similar} on $\sii(K)$ to
the equivalent problem
\begin{equation}\label{heat.similar.HO}
  v_s
  + \big(\!-i\nabla_{\!y}-A_s(y)\big)^2 v
  + \frac{\ |y|^2}{16} \, v
  = 0
\end{equation}
on $\sii$. We observe the presence of the harmonic-oscillator potential
that makes the spectrum of the generator to be discrete.
Furthermore, recalling the spectrum of 
the quantum-harmonic-oscillator Hamiltonian
(see any textbook on quantum mechanics, \eg, \cite[Sec.~2.3]{Griffiths})
and the relationship between~$s$ and~$t$ from~\eqref{SST},
it is clear now why we get the~$1/2$ for the decay rate
in Theorem~\ref{Thm.rate} if $A=0$.
It remains to analyse the effect of the presence of the magnetic field.
\smallskip \\
\textbf{IV (taking $s \to \infty$).}
We look at the asymptotic behaviour of~\eqref{heat.similar.HO}
as the self-similar time~$s$ tends to infinity.
The scaling coming from the self-similarity transformation
is such that
\begin{equation}\label{AB-limit}
  B_s := \rot A_s
  \xrightarrow[s \to \infty]{}
  2\pi \, \Phi_B \, \delta =: B_\infty
\end{equation}
in the sense of distribution, where~$\delta$ is the Dirac delta function
and~$\Phi_B$ is total flux~\eqref{flux}.
The right hand side of~\eqref{AB-limit} is the Aharonov-Bohm field
to which there corresponds an explicit vector potential~$A_\infty$
having a singularity at~$0$.
Taking this limit into account,
it is expectable that~\eqref{heat.similar.HO}
will be approximated for large times by the Aharonov-Bohm problem
\begin{equation}\label{heat.similar.HO.limit}
  w_s
  + \big(\!-i\nabla_{\!y}-A_\infty(y)\big)^2 w
  + \frac{\ |y|^2}{16} \, w
  = 0
  \,.
\end{equation}
This evolution equation is explicitly solvable
in terms of special functions.
In particular, it is clear that the solutions
decay exponentially with the decay rate given by
the lowest eigenvalue of the (time-independent) generator
which explicitly reads
\begin{equation}\label{lev}
  (1+\beta)/2
  \,.
\end{equation}
It remains to observe that the exponential decay rate transfers
to a polynomial one after returning to the original time~$t$
and that the number~\eqref{lev}
is just that of the bound of Theorem~\ref{Thm.rate}
in the presence of magnetic field.

\smallskip
Two comments are in order.
First, although it is very natural to expect,
we do not establish any theorem
that solutions of~\eqref{heat.similar.HO} can be approximated
by those of~\eqref{heat.similar.HO.limit} as $s \to \infty$.
We only show a strong-resolvent convergence for operators
related to their generators (Proposition~\ref{Prop.strong}).
This is, however, sufficient to prove Theorem~\ref{Thm.rate}
with help of energy estimates.
Proposition~\ref{Prop.strong} is probably the most
significant auxiliary result of the paper
and we believe it is interesting in its own right.
Second, in the proof of Proposition~\ref{Prop.strong}
we heavily rely on the existence of
the Hardy-type inequality~\eqref{Hardy}.

\subsection{The content of the paper}
The organization of this paper is as follows.

In the following Section~\ref{Sec.Pre}
we give a precise definition of the magnetic 
Schr\"o\-dinger operator~\eqref{Hamiltonian}
and the associated semigroup.
We also briefly discuss there the case 
of zero magnetic field,
obtaining, \emph{inter alia},
the first statement of Theorem~\ref{Thm.rate}
(\cf~Proposition~\ref{Prop.straight} and Corollary~\ref{Corol.straight}).

The main body of the paper is represented by Section~\ref{Sec.ss}
where we develop the method of self-similar solutions
to get the improved decay rate of Theorem~\ref{Thm.rate}
as described above.
Here we point out Section~\ref{Sec.ess} where a convergence
of operators associated with the regular and Aharonov-Bohm
fields of~\eqref{AB-limit} is proved (\cf~Proposition~\ref{Prop.strong}).
Finally, in Section~\ref{Sec.alt} we establish
an alternative version of Theorem~\ref{Thm.rate}
which provides a global (in time) upper bound to
the heat semigroup (\cf~Theorem~\ref{Thm.rate.alt}).

The paper is concluded in Section~\ref{Sec.end}
by referring to some open problems.

\section{Preliminaries}\label{Sec.Pre}
%
Let the local magnetic field~$B$ be represented
by a real-valued continuous function of compact support in~$\Real^2$.
It is well known that the corresponding magnetic vector potential~$A$
is not uniquely determined by the formula
\begin{equation}\label{rot}
  B = \rot A \equiv \partial_1 A_2 - \partial_2 A_1
  \,.
\end{equation}
This freedom, reflecting the gauge invariance of the physical theory,
enables us to work with the particular choice
\begin{equation}\label{potential}
  A(x) := (-x_2,x_1) \int_0^1 B(\tau x_1,\tau x_2) \, \tau \, d\tau
\end{equation}
without loss of any generality.
It corresponds to the so-called \emph{multipolar gauge} in molecular electrodynamics,
and it is also known as the \emph{transversal gauge} because
$
  x \cdot A(x) = 0
$.
For smooth~$B$ it is straightforward to verify the validity of~\eqref{rot}.
If~$B$ is merely continuous, the identity~\eqref{rot}
can be checked in the sense of distributions
by approximating~$B$ by a sequence of mollifiers.

In polar coordinates $(r,\theta) \in [0,\infty) \times [0,2\pi)$,
related to the Cartesian system by
$
  x = (r\cos\theta,r\sin\theta)
$,
the magnetic potential~\eqref{potential}
acquires a more transparent form
\begin{equation}\label{potential.polar}
  A(x) = (-\sin\theta,\cos\theta) \ \frac{\alpha(r,\theta)}{r}
\end{equation}
with the notation
\begin{equation*}
  \alpha(r,\theta) :=
  \int_0^r B(\tau\cos\theta,\tau\sin\theta) \, \tau \, d\tau
\end{equation*}
Since~$B$ is supposed to be compactly supported,
the function~$\alpha$ is bounded and~$A$ is vanishing at infinity.
Moreover, the limit 
\begin{equation}\label{alpha.infinity}
  \alpha_\infty(\theta) := \lim_{r\to\infty} \alpha(r,\theta)
\end{equation}
is well defined for every $\theta\in[0,2\pi)$.

The quantity
$$
  \Phi(r) := \frac{1}{2\pi} \, \int_0^{2\pi }\alpha(r,\theta) \, d\theta
$$
has the physical meaning of the \emph{magnetic flux}
through the disc of radius~$r$ centred at the origin of~$\Real^2$.
Since the field~$B$ is zero outside such a disc of sufficiently large radius, 
$\Phi(r)=\Phi_B$ for all~$r$ large enough,
where~$\Phi_B$ is the \emph{total magnetic flux} introduced in~\eqref{flux}.
As pointed out in the introduction, the set of integers correspond
to \emph{flux quanta} for our choice of physical constants
(\cf~\cite[Sec.~10.2.4]{Griffiths}). 

The magnetic Hamiltonian~\eqref{Hamiltonian} is defined as
the self-adjoint operator associated on $\sii:=\sii(\Real^2)$
with the closed quadratic form
\begin{equation}\label{form}
  h_B[\psi] := \|\nabla\psi-iA\psi\|^2
  \,, \qquad
  \psi \in \Dom(h_B) := W^{1,2}
  \,,
\end{equation}
where~$\|\cdot\|$ denotes the norm in~$\sii$
and $W^{1,2}:=W^{1,2}(\Real^2)$ is the usual Sobolev space.
Using the transverse gauge and polar coordinates,
we can write
\begin{equation}\label{form.polar}
  h_B[\psi] = \|\partial_r\psi\|^2
  + \|r^{-1}(\partial_\theta\psi-i\alpha\psi)\|^2
  \,.
\end{equation}
\begin{Remark}
By the gauge invariance mentioned above,
it makes sense to distinguish the Hamiltonian
by the subscript~$B$ rather than~$A$.
Indeed, all the Hamiltonians corresponding to different gauges
are unitarily equivalent
(\cf~\cite[Prob.~4.53]{Griffiths}).
\end{Remark}

Since $H_B-H_0$ is a relatively form compact perturbation
\cite[Prob.~XIII.39]{RS4}
of the free Hamiltonian~$H_0$,
we have
\begin{equation}\label{spectrum}
  \sigma(H_B) = [0,\infty)
  \,.
\end{equation}

As usual, we consider the weak formulation of
the parabolic problem~\eqref{heat.intro}
subjected to the initial datum $u_0\in\sii$.
We say a Hilbert space-valued function
$
  u \in \sii_\mathrm{loc}\big((0,\infty);W^{1,2}\big)
$,
with the weak derivative
$
  u' \in \sii_\mathrm{loc}\big((0,\infty);W^{-1,2}\big)
$,
is a (global) solution of~\eqref{heat.intro} provided that
\begin{equation}\label{heat.weak}
  \big\langle v,u'(t)\big\rangle + h_B\big(v,u(t)\big) = 0
\end{equation}
for each $v \in W^{1,2}$ and a.e.\ $t\in[0,\infty)$,
and $u(0)=u_0$.
Here $h_B(\cdot,\cdot)$ denotes the sesquilinear form
associated with~\eqref{form}
and $\langle\cdot,\cdot\rangle$
stands for the pairing of $W^{1,2}$ and its dual $W^{-1,2}$.
With an abuse of notation, we denote by the same symbol~$u$
both the function on $\Real^2\times(0,\infty)$
and the mapping $(0,\infty) \to W^{1,2}$.

Standard semigroup theory implies that there indeed exists
a unique solution of~\eqref{heat.weak} that belongs to
$C^0\big([0,\infty),\sii\big)$.
More precisely, the solution is given by $u(t) = e^{-H_B t} u_0$,
where~$e^{-H_B t}$ is the semigroup associated with~$H_B$.
It is important to stress that the semigroup is not positivity-preserving,
which in particular means that we have to work with
complex functional spaces when studying~\eqref{heat.weak}.

The spectral mapping theorem together with~\eqref{spectrum} yields
$
  \|e^{-H_B t}\|_{\sii\to\sii} = 1
$
for each time $t \geq 0$.
To observe an additional large-time decay of the solutions~$u(t)$,
we have to restrict the class of initial data~$u_0$ to a subspace of~$\sii$.
In this paper we consider the weighted space $\sii(K):=\sii(\Real^2,K(x)\,dx)$
with the weight given by~\eqref{weight}.
The norm in $\sii(K)$ will be denoted by $\|\cdot\|_K$.

The case of the semigroup of the free Hamiltonian~$H_0$
is well understood because of the explicit knowledge
of its heat kernel:
$$
  p_0(x,x',t) := (4 \pi t)^{-1} e^{-|x-x'|/(4t)}
  \,.
$$
Using this expression,
it is straightforward to establish the following bounds:
\begin{Proposition}\label{Prop.straight}
There exists a constant~$C$ such that for every $t \geq 1$,
$$
  C^{-1} \, t^{-1/2}
  \leq
  \|e^{-H_0 t}\|_{\sii(K) \to \sii}
  \leq
  C \, t^{-1/2}
  \,.
$$
\end{Proposition}
\begin{proof}
Using the Schwarz inequality, we get
$$
  \|e^{-H_0 t} u_0\|^2 \leq \|u_0\|_K^2
  \int_{\Real^2\times\Real^2}
  \frac{|p_0(x,x',t)|^2}{K(x')} \, dx \, dx'
$$
for every $u_0 \in \sii(K)$,
where the integral can be evaluated explicitly
and yields the upper bound.
To get the lower bound, it is enough to find
one particular $u_0\in\sii(K)$ such that
$\|e^{-H_0 t} u_0\| \geq C^{-1} \, t^{-1/2} \, \|u_0\|_K $.
It is just a matter of explicit computations to check
that it is the case for the choice $u_0 = K^{-\beta}$ with any $\beta>1/2$.
\end{proof}

As a consequence of this proposition,
we get:
\begin{Corollary}\label{Corol.straight}
We have $\gamma_0=1/2$.
\end{Corollary}

This establishes the trivial part of Theorem~\ref{Thm.rate}.
On the other hand, no explicit formula for the heat kernel
is available if $B\not=0$ and a more advanced technique
is needed to study the decay rate~$\gamma_B$.

\section{The self-similarity transformation}\label{Sec.ss}
%
Our method to study the asymptotic behaviour of
the heat equation~\eqref{heat.intro} in the presence
of magnetic field is to adapt
the technique of self-similar solutions
used in the case of the heat equation in the whole Euclidean space
by Escobedo and Kavian~\cite{Escobedo-Kavian_1987}
to the present problem.
We closely follow the approach of the recent papers~\cite{KZ1,KZ2},
where the technique is applied to twisted waveguides
in three and two dimensions, respectively.

\subsection{An equivalent time-dependent problem}
We consider a unitary transformation~$U$ on~$\sii$
which associates to every solution
$
  u \in \sii_\mathrm{loc}\big((0,\infty),dt;\sii\big)
$
of~\eqref{heat.weak}
a function~$\tilde{u}=U u$
in a new $s$-time weighted space
$
  \sii_\mathrm{loc}\big((0,\infty),e^s ds;\sii\big)
$
via~\eqref{SST}.
The inverse change of variables is given by
\begin{equation}\label{SST.inverse}
  u(x,t)
  = (t+1)^{-1/2} \, \tilde{u}\big((t+1)^{-1/2}x,\log(t+1)\big)
  \,.
\end{equation}
Note that the original space-time variables $(x,t)$
are related to the ``self-similar'' space-time	 variables $(y,s)$ 
via the relations
\begin{equation}\label{space-times}
  (x,t) = (e^{s/2}y,e^s-1)
  \,, \qquad
  (y,s) = \big((t+1)^{-1/2}x,\log(t+1)\big)
  \,.
\end{equation}
Hereafter we consistently use the notation for respective variables
to distinguish the two space-times. 

It is easy to check that this change of variables leads
to the evolution equation~\eqref{heat.similar},
subject to the same initial condition as~$u$ in~\eqref{heat.intro}.
More precisely, the weak formulation~\eqref{heat.weak}
transfers to
\begin{equation}\label{heat.weak.similar}
  \big\langle
  \tilde{v}, \tilde{u}'(s)
  -\mbox{$\frac{1}{2}$} \, y \cdot \nabla_{\!y}\tilde{u}(s)
  \big\rangle
  + \tilde{a}_{s}\big(\tilde{v},\tilde{u}(s)\big) = 0
\end{equation}
for each $\tilde{v} \in W^{1,2}$ and a.e.~$s\in[0,\infty)$,
with $\tilde{u}(0) = \tilde{u}_0 := U u_0 = u_0$.
Here~$\tilde{a}_{s}(\cdot,\cdot)$ denotes the sesquilinear form
associated with
\begin{equation*}
  \tilde{a}_{s}[\tilde{u}] :=
  \|\nabla\tilde{u}-i A_s\tilde{u}\|^2
  - \frac{1}{2} \, \|\tilde{u}\|^2
  \,, \qquad
  \tilde{u} \in \Dom(\tilde{Q}_{s}) := W^{1,2}
  \,,
\end{equation*}
where~$A_s$ has been introduced in~\eqref{sigma}.

\begin{Remark}
Note that~\eqref{heat.similar} (or \eqref{heat.weak.similar}) 
is a parabolic equation with $s$-time-dependent coefficients.
The same occurs and has been previously analysed
in the twisted waveguides \cite{KZ1,KZ2}
and also for a convection-diffusion equation in the whole space
but with a variable diffusion coefficient
\cite{Escobedo-Zuazua_1991,Duro-Zuazua_1999}.
A careful analysis of the behaviour of the underlying elliptic operators
as~$s$ tends to infinity leads to a sharp decay rate for its solutions.
\end{Remark}

Since~$U$ acts as a unitary transformation on~$\sii$,
it preserves the space norm of solutions
of~\eqref{heat.intro} and~\eqref{heat.similar}, \ie,
\begin{equation}\label{preserve}
  \|u(t)\|=\|\tilde{u}(s)\|
  \,.
\end{equation}
This means that we can analyse the asymptotic time behaviour
of the former by studying the latter.

However, the natural space to study the evolution~\eqref{heat.similar}
is not~$\sii$ but rather the weighted space $\sii(K)$.
Following the approach of~\cite{KZ1}
based on a theorem of J.~L.~Lions~\cite[Thm.~X.9]{Brezis_FR}
about weak solutions of parabolic equations
with time-dependent coefficients,
it can be shown that~\eqref{heat.similar}
is well posed in the scale of Hilbert spaces
\begin{equation}\label{scale}
  W^{1,2}(K)
  \subset \sii(K) \subset
  W^{1,2}(K)^*
  \,,
\end{equation}
with
$$
  W^{1,2}(K)
  := \left\{ \tilde{u}\in \sii(K) \ | \
  \nabla\tilde{u}\in \sii(K)
  \right\}
  .
$$

More precisely, choosing $\tilde{v}:= K v$
for the test function in~\eqref{heat.weak.similar},
where $v \in C_0^\infty(\Omega_0)$ is arbitrary,
we can formally cast~\eqref{heat.weak.similar}
into the form
\begin{equation}\label{heat.weak.weighted}
  \big\langle v, \tilde{u}'(s) \big\rangle_K
  + a_s\big(v,\tilde{u}(s)\big) = 0
  \,.
\end{equation}
Here $\langle\cdot,\cdot\rangle_K$
denotes the pairing of $W^{1,2}(K)$ and $W^{1,2}(K)^*$,
and
\begin{equation}\label{form.a}
  a_s(v,\tilde{u}) :=
  \big(\nabla v - i A_s v,
  \nabla \tilde{u} - i A_s \tilde{u}
  \big)_K
  - \frac{1}{2} \, \big(y \;\! v,
  i A_s \tilde{u}\big)_K
  - \frac{1}{2} \, \big(v,\tilde{u}\big)_K
  \,,
\end{equation}
with $(\cdot,\cdot)_K$ denoting the inner product in $\sii(K)$.
Note that~$a_s$ is not a symmetric form.

By ``formally'' we mean that the formulae are meaningless in general,
because the solution~$\tilde{u}(s)$ and its derivative~$\tilde{u}'(s)$
may not belong to $W^{1,2}(K)$ and $W^{1,2}(K)^*$, respectively.
The justification of~\eqref{heat.similar}
being well posed in the scale~\eqref{scale}
consists basically in checking the boundedness
and a coercivity of the form~$a_s$ defined on $\Dom(a_s):=W^{1,2}(K)$.
It is straightforward by noticing that~$A_s$ is bounded
for each~$s$ fixed.
We refer to \cite[Sec.~5.4]{KZ1} for more details.

\subsection{Reduction to a spectral problem}
Choosing $v := \tilde{u}(s)$ in~\eqref{heat.weak.weighted}
and combining the obtained equation with its conjugate version,
we arrive at the identity
\begin{equation}\label{formal}
  \frac{1}{2} \frac{d}{ds} \|\tilde{u}(s)\|_{K}^2
  = - \hat{l}_s[\tilde{u}(s)]
  \,,
\end{equation}
where $\hat{l}_s[\tilde{u}] := \Re\{a_s[\tilde{u}]\}$,
$\tilde{u} \in \Dom(\hat{l}_s) := \Dom(a_s) = W^{1,2}(K)$
(independent of~$s$).
Observing that the real part of the non-symmetric term
in~\eqref{form.a} vanishes, we get
\begin{equation*}
  \hat{l}_s[\tilde{u}] :=
  \big\|\nabla \tilde{u} - i A_s \tilde{u}\big\|_K
  - \frac{1}{2} \, \|\tilde{u}\|_{K}^2
  \,.
\end{equation*}
It remains to analyse the coercivity of~$\hat{l}_s$.

More precisely, as usual for energy estimates,
we replace the right hand side of~\eqref{formal}
by the spectral bound, valid for each fixed $s \in [0,\infty)$,
\begin{equation}\label{spectral.reduction}
  \forall \tilde{u} \in \Dom(\hat{l}_s) \;\!, \qquad
  \hat{l}_s[\tilde{u}]
  \geq \lambda(s) \, \|\tilde{u}\|_{K}^2
  \,,
\end{equation}
where~$\lambda(s)$ denotes the lowest point in the spectrum of
the self-adjoint operator~$\hat{L}_s$
associated on~$\sii(K)$ with~$\hat{l}_s$.
Then~\eqref{formal} together with~\eqref{spectral.reduction} implies
the exponential bound
\begin{equation}\label{spectral.reduction.integral}
  \forall s \in [0,\infty) \;\!, \qquad
  \|\tilde{u}(s)\|_{K}
  \leq \|\tilde{u}_0\|_{K} \
  e^{-\int_0^s \lambda(\tau) \, d\tau}
  \,.
\end{equation}
Finally, recall that the exponential bound in~$s$
transfers to a polynomial bound in the original time~$t$,
\cf~\eqref{space-times}.
In this way, the problem is reduced to a spectral analysis
of the family of operators $\{\hat{L}_s\}_{s \geq 0}$.

\subsection{Study of the spectral problem}
In order to investigate the operator~$\hat{L}_s$ in~$\sii(K)$,
we first map it into a unitarily equivalent operator
$L_s := \mathcal{U} \hat{L}_s \mathcal{U}^{-1}$
on~$\sii$ via the unitary transform
$$
  \mathcal{U}\;\!\tilde{u} := K^{1/2}\,\tilde{u}
  \,.
$$
By definition, $L_s$~is the self-adjoint operator
associated on~$\sii$ with the quadratic form
$
  l_s[v] := \hat{l}_s[\mathcal{U}^{-1}v]
$,
$
  v \in \Dom(l_s) := \mathcal{U}\,\Dom(\hat{l}_s)
$.
A straightforward calculation yields
\begin{equation}\label{J0.form}
  l_s[v]
  = \|\nabla v - i A_s v\|^2
  + \frac{1}{16} \, \|y v\|^2
  \,, \qquad
  v \in \Dom(l_s)
  = W^{1,2}(\Real^2) \cap \sii\big(\Real^2, |y|^2 \, dy\big)
  \,.
\end{equation}
In particular, $\Dom(l_s)$~is independent of~$s$.
In the distributional sense, we can write
\begin{equation}\label{Hamiltonian.ss}
  L_s =
  (-i\nabla - A_s)^2 + \frac{1}{16} \, |y|^2
  \,.
\end{equation}

For the trivial magnetic field $B=0$, the operator~$L_s$ coincides with
the Hamiltonian of the quantum harmonic oscillator
\begin{equation}\label{HO}
  L_\mathrm{HO} :=
  -\Delta + \frac{1}{16} \, |y|^2
  \qquad \mbox{in} \qquad
  \sii(\Real^2)
\end{equation}
(\ie\ the Friedrichs extension
of this operator initially defined on $C_0^\infty(\Real^2)$).
In any case, the form domain of $L_\mathrm{HO}$ coincides with~$\Dom(l_s)$.
Recall the well known fact that~$L_\mathrm{HO}$ has compact resolvent
(see, \eg, \cite[Thm.~XIII.67]{RS4}) and thus purely discrete spectrum.
We conclude that, for any field~$B$,
$L_s$~has compact resolvent, too.
In particular, $\lambda(s)$~represents the lowest
eigenvalue of~$L_s$.

By the diamagnetic inequality \cite[Thm.~2.1.1]{Fournais-Helffer_2009}
\begin{equation}\label{diamagnetic}
  \big|(-i\nabla - A_s) v\big| \geq \big|\nabla|v|\big|
\end{equation}
valid for every $v \in \sii$ almost everywhere
and by the minimax principle, it follows that~$\lambda(s)$
is bounded from below by the lowest eigenvalue of~$L_\mathrm{HO}$.
Using the explicit knowledge of the spectrum of~$L_\mathrm{HO}$
(see, \eg, \cite[Sec.~2.3]{Griffiths}),
we get
\begin{equation}\label{diamagnetic.bound}
  \forall s\in[0,\infty)
  \,, \qquad
  \lambda(s) \geq 1/2
  \,.
\end{equation}
Obviously, there is an equality if $B=0$.
Moreover, the bound becomes sharp in the limit $s\to\infty$
for magnetic fields with the total flux being an integer.
\begin{Proposition}\label{Prop.sharp}
Let $\Phi_B \in \Int$. Then
$$
  \lambda(\infty) := \lim_{s\to\infty} \lambda(s) = 1/2
  \,.
$$
\end{Proposition}
\begin{proof}
In view of~\eqref{diamagnetic.bound}, it is enough to establish
an upper bound to~$\lambda(s)$ that goes to~$1/2$ as $s\to\infty$.  
We obtain it by constructing a suitable trial function
in the variational characterization of the eigenvalue.	

For every $n\in\Nat$ larger than~$1$, 
let us define the function
$$
\displaystyle
\eta_n(r) :=
\begin{cases}
  1 
  & \mbox{if}\quad r>1/n \,,
  \\
  \displaystyle
  \frac{\log(n^2 r)}{\log n}
  & \mbox{if}\quad r\in[1/n^2,1/n] \,,
  \\
  0
  & \mbox{otherwise} \,.
\end{cases}
$$
We have $|\eta_n|\leq1$ for all $n \geq 2$, 
$\eta_n(r)$~converges as $n\to\infty$
to~$1$ for every $r\in(0,\infty)$ and 
\begin{equation}\label{eta.convergence}
  \|\eta_n'\|_{\sii((0,\infty),r dr)}^2
  = (\log n)^{-1}
  \xrightarrow[n\to\infty]{}
  0
  \,.
\end{equation}
Hence, the functions $(r,\theta) \mapsto \eta_n(r)$
in polar coordinates in~$\Real^2$ represent
a convenient approximation of the constant function~$1$ 
in $W_\mathrm{loc}^{1,2}(\Real^2)$ 
when a singularity in the origin has to be avoided. 

In polar coordinates in~$\Real^2$, we then introduce 
$$
  \psi_{n}(r,\theta) := \varphi(r) \, \eta_n(r) \, 
  e^{i\int_0^\theta \alpha_\infty(\tau)\,d\tau}
  \,,
$$
where $(r,\theta) \mapsto \varphi(r) := e^{-r^2/8}$ is an eigenfunction 
of~$L_\mathrm{HO}$ corresponding to its lowest eigenvalue~$1/2$
(\cf~Proposition~\ref{Prop.AB-spectrum}) 
and~$\alpha_\infty$ is defined in~\eqref{alpha.infinity}.
Since~$\psi_n$ vanishes in the vicinity of the origin 
and decays sufficiently fast at the infinity,
we clearly have $\psi_n \in \Dom(l_s)$ 
for all $n \geq 2$ and $s \geq 0$.
The dominated convergence theorem implies
$$
  \|\psi_n\|^2 \xrightarrow[n\to\infty]{} 
  2\pi \int_0^\infty \varphi(r)^2 \, r \, dr
  = 4\pi 
  \,.
$$
At the same time,
a straightforward calculation yields
\begin{align*}
  l_s[\psi_n] - \mbox{$\frac{1}{2}$} \|\psi_n\|^2
  \ = \ & \int_0^\infty \int_0^{2\pi}
  \varphi(r)^2 \, \eta_n(r)^2 \,
  \big|\alpha_\infty(\theta)-\alpha(e^{s/2} r,\theta)\big|^2 \,
  d\theta \, \frac{dr}{r}
  \\
  & + 2\pi \int_0^\infty \varphi(r)^2 \, \eta_n'(r)^2 \, r \, dr
  \,.
\end{align*}
Since $\alpha(e^{s/2} r,\theta)$ converges as $s\to\infty$
to~$\alpha_\infty$ for every $(r,\theta)\in(0,\infty)\times[0,2\pi)$  
and~$\alpha$ is a bounded function, 
the double integral tends to zero as $s\to\infty$
for every $n \geq 2$ fixed due to the dominated convergence theorem.  
The one-dimensional integral is independent of~$s$
and it goes to zero as $n\to\infty$ 
due to~\eqref{eta.convergence} and $|\varphi|<1$.  

By~\eqref{diamagnetic.bound} 
and by the Rayleigh-Ritz variational formula for~$\lambda(s)$,
we have the bounds
$$
  0 
  \leq \lambda(s) - \frac{1}{2} \leq 
  \frac{l_s[\psi_n]}{\|\psi_n\|}
$$
for every $s \geq 0$.
Fixing $n \geq 2$ and taking the limit $s\to\infty$, 
we therefore get 
$$
  0 
  \leq \lambda(\infty) - \frac{1}{2} \leq
  \frac{\int_0^\infty \varphi(r)^2 \, \eta_n'(r)^2 \, r \, dr}
  {\int_0^\infty \varphi(r)^2 \, r \, dr}
  \,.
$$ 
Since the right hand side provides an upper bound 
which can be made arbitrarily small by choosing~$n$ 
sufficiently large, we conclude with the statement of the proposition. 
\end{proof} 

Proposition~\ref{Prop.sharp} explains 
why the estimate~\eqref{spectral.reduction.integral}
does not lead to any improved polynomial decay rate if $\Phi_B \in \Int$. 
To get the improved decay rate for $\Phi_B \not\in \Int$,
we need to show that a better inequality than~\eqref{diamagnetic.bound}
is available, at least asymptotically as $s\to\infty$.

\subsection{The asymptotic behaviour of the spectrum}\label{Sec.ess}
This subsection is devoted to a study of the asymptotic
behaviour of the operator~$L_s$ as $s\to\infty$.
Since the limit directly leads to the improved decay rate,
it is probably the most essential part of the paper. 

We recall the expression~\eqref{potential.polar} 
for the magnetic potential~$A$ in the transverse gauge
and conventionally expressed in polar coordinates.
It is clear that the scaling~\eqref{sigma} acts in such a way that
\begin{equation}\label{AB-limit-potential}
  A_s \xrightarrow[s\to\infty]{} A_\infty
  \qquad
  \mbox{locally uniformly in $\Real^2\setminus\{0\}$}
  \,,
\end{equation}
where (recall~\eqref{alpha.infinity})
\begin{equation}\label{AB-potential}	
  A_\infty(y) 
  := (-\sin\theta,\cos\theta) \ \frac{\alpha_\infty(\theta)}{r}
\end{equation}
is well defined off the origin.
Here, to save symbols, we use the same notation 
$(r,\theta)\in[0,\infty)\times[0,2\pi)$ 
for polar coordinates in the ``self-similar'' $y$-plane as well.  
Contrary to~$A_s$, which is a bounded function for every $s\in[0,\infty)$,  
$A_\infty$~possesses a singularity at the origin.  
As already revealed in the introduction, 
$A_\infty$~represents the magnetic potential in the transverse gauge
for the Aharonov-Bohm field~$B_\infty$ defined in~\eqref{AB-limit}.
Indeed, $\rot A_\infty = 2\pi \Phi_B \, \delta$ in the sense of distributions.

\subsubsection{The Aharonov-Bohm Hamiltonian}
In view of~\eqref{AB-limit-potential},
it is reasonable to expect that~$L_s$ 
converges in a suitable sense as $s\to\infty$ 
to the operator  
\begin{equation}\label{Hamiltonian.AB}
  L_\infty :=
  (-i\nabla - A_\infty)^2 + \frac{1}{16} \, |y|^2
  \,.
\end{equation}
However, one has to be careful when giving a meaning to~$L_\infty$
because of the singularity of~$A_\infty$. 
Indeed, already for the operator $(-i\nabla - A_\infty)^2$
with a non-trivial rotationally-symmetric Aharonov-Bohm field,
\ie~$\alpha_\infty(\theta)=\const\not\in\Int$, it is well known 
that the symmetric operator defined by this differential 
expression on the minimal domain 
$C_0^\infty(\Real^2\setminus\{0\})$ is not essentially self-adjoint
and that there are infinitely many self-adjoint extensions determined 
by imposing boundary conditions at the origin
(\cf~\cite{Pankrashkin-Richard} 
for a modern treatment in terms of boundary triplets).

It turns out that the ``limit'' of~$L_s$ as $s\to\infty$ corresponds
to the Friedrichs extension of~$L_\infty$, 
\ie, we understand~\eqref{Hamiltonian.AB}
as the self-adjoint operator associated on~$\sii$ 
with the quadratic form 
\begin{equation}\label{J0.form.AB}
  l_\infty[v]
  := \|\nabla v - i A_\infty v\|^2
  + \frac{1}{16} \, \|y v\|^2
  \,, \qquad
  v \in \Dom(l_\infty)
  := \overline{C_0^\infty(\Real^2\setminus\{0\})}^{\|\cdot\|_{l_\infty}}
  \,.
\end{equation}
Here the norm $\|\cdot\|_{l_\infty}$ with respect to which 
the closure is taken is defined by
%
$  
  \|v\|_{l_\infty}^2 := l_\infty[v]+\|v\|^2 
$.
%
We call~$L_\infty$ simply the Aharonov-Bohm Hamiltonian,
although it contains the additional
quantum-harmonic-oscillator potential term
with respect to the classical Aharonov-Bohm Hamiltonian.

Now we specify the domain of the form~$l_\infty$.
Observing the structure of~$l_\infty$ in polar coordinates 
(\cf~\eqref{form.polar}) 
and being inspired by~\cite{Laptev-Weidl_1999}, 
we consider the ordinary differential self-adjoint operator~$K$
on the Hilbert space $\sii((0,2\pi))$ defined by
$$
  K\phi := i\phi' + \alpha_\infty \phi 
  \,, \qquad
  \phi\in\Dom(K) := \big\{\phi\in W^{1,2}((0,2\pi))\,|\,\phi(0)=\phi(2\pi)\big\}
  \,.
$$
The spectral problem for~$K$ is explicitly solvable.
In particular, we have 
$
  \sigma(K) = \{m + \Phi_B \}_{m\in\Int}
$  
and the corresponding set of (normalized) eigenfunctions is given by
$$
  \phi_m(\theta) := 
  (2\pi)^{-1/2} \
  e^{
  -i
  \left[(m + \Phi_B)\theta-\int_0^\theta \alpha_\infty(\tau) \, d\tau\right]
  }
  \,.
$$ 
Decomposing~$L_\infty$ into this complete orthonormal set of $\sii((0,2\pi))$,
we thus get 
\begin{equation}\label{decomposition}
  L_\infty 
  \ \simeq \ \bigoplus_{m\in\Int} 
  \left(
  - \frac{1}{r}\frac{d}{dr} r \frac{d}{dr} 
  + \frac{(m + \Phi_B)^2}{r^2} + \frac{r^2}{16}
  \right)
  \,.
\end{equation}
Here the symbol~$\simeq$ stands for the unitary-equivalence relation,
the ordinary differential operators on the right hand side 
should be understood as the Friedrichs extensions in $\sii((0,\infty),r\,dr)$ 
of these operators defined initially on $C_0^\infty((0,\infty))$
and considered as acting on the respective space
$\sii((0,\infty),r\,dr) \otimes \{\phi_m\}$ by natural isomorphisms.

It follows that~$L_\infty$ is unitarily equivalent 
to the harmonic-oscillator Hamiltonian~$L_\mathrm{HO}$ 
from~\eqref{HO} if $\Phi_B\in\Int$
(this can be also understood in terms of a gauge invariance).  
Recall that the form domain of~$L_\mathrm{HO}$ 
coincides with $\Dom(l_s)$ given in~\eqref{J0.form}.

On the other hand, the decomposition~\eqref{decomposition} immediately yields 
the Hardy-type inequality
\begin{equation}\label{Hardy.AB}
  L_\infty 
  \geq 
  \frac{\beta^2}{\,|\cdot|^2} 
  \qquad \mbox{with} \qquad
  \beta:=\dist(\Phi_B,\Int)
\end{equation}
in the sense of quadratic forms,
which becomes non-trivial whenever $\Phi_B\not\in\Int$. 
This inequality can be also obtained as a consequence 
of the Laptev-Weidl Hardy-type inequality 
\cite[Thm.~3]{Laptev-Weidl_1999},
representing an elegant modification of~\eqref{Hardy} 
for the Aharonov-Bohm field.

Here we use~\eqref{Hardy.AB} just to conclude with
\begin{equation}\label{topology}
    \Dom(l_\infty) = 
    W^{1,2}(\Real^2) \cap \sii\big(\Real^2, |y|^2 \, dy\big)
    \cap \sii\big(\Real^2, |y|^{-2} \, dy\big)
    \quad\mbox{if} \quad 
    \Phi_B\not\in\Int
    \,.
\end{equation}
Indeed, if $\Phi_B\not\in\Int$,
the inequality~\eqref{Hardy.AB} enables one
to handle the mixed terms after developing 
the second norm in~\eqref{form.polar} and to show 
that the norm~$\|\cdot\|_{l_\infty}$ is actually equivalent
to that induced by the form
$
  \|\nabla v\|^2 + \|y v\|^2 + \||y|^{-1}v\|^2
$.
It follows that the elements of~$\Dom(l_\infty)$ 
have to vanish at the origin if $\Phi_B\not\in\Int$
(while for $\Phi_B\in\Int$ they can even diverge there).

In any case, the form domain~$\Dom(l_\infty)$ is compactly embedded in~$\sii$,
so that the operator~$L_\infty$ has compact resolvent and purely discrete spectrum.
The latter is explicitly computable as we show now.
\begin{Proposition}\label{Prop.AB-spectrum}
The spectrum of~$L_\infty$ is given by the discrete set
$$
  \sigma(L_\infty) = 
  \left\{
  n + \frac{1+|m + \Phi_B|}{2}
  \right\}_{\! n\in\Nat, \, m\in\Int}
  \,.
$$
The corresponding set of (unnormalized) eigenfunctions is given 
in polar coordinates by
\begin{equation}\label{Laguerre}
  r^{|m + \Phi_B|} \, e^{-r^2/8} \, \mathcal{L}_n^{|m + \Phi_B|}(r^2/4) 
  \ \phi_m(\theta)
  \,,
\end{equation}
where $\mathcal{L}_n^{\mu}$ denotes the generalized Laguerre polynomial.
\end{Proposition}
\begin{proof}
As a consequence of the decomposition~\eqref{decomposition}, we get 
$$
  \sigma(L_\infty) = \bigcup_{m\in\Int} \, \sigma(L_\infty^m)
  \,,
$$
where~$L_\infty^m$ denotes the ordinary-differential operator
on the right hand side of~\eqref{decomposition}.
Consequently, the problem reduces to the study of 
the spectrum of each~$L_\infty^m$.
Fix $m \in \Int$.

The eigenvalue problem $L_\infty^m\psi=E\psi$ 
leads to an ordinary-differential equation that 
can be cast into Whittaker's equation \cite[13.1.31]{AS}.
The latter admits solutions in terms of 
confluent hypergeometric functions~$M$ and~$U$, \cf~\cite[Sec.~13.1]{AS}. 
Checking the asymptotic behaviours of the functions 
around the origin and at the infinity, 
it is straightforward to verify that none of the solutions
belongs to~$\Dom(L_\infty^m)$ because of singularities,
unless the following quantization condition 
for the eigenvalues
$$
  E_n = n + \frac{1+|m + \Phi_B|}{2}
  \,, \qquad
  n \in \Nat \equiv \{0,1,2,\dots\} \,,
$$
holds true. 
Under this condition, the confluent hypergeometric functions reduce 
to generalized Laguerre's polynomials \cite[Sec.~22]{AS} 
and the latter lead to a complete orthonormal system of 
eigenfunctions for each~$L_\infty^m$.

More straightforwardly, it follows from \cite[22.6.18]{AS} 
that~\eqref{Laguerre} is a solution of
$L_\infty^m\psi=E_n\psi$ for each fixed $m \in \Int$.
Then it remains to recall the classical fact 
(see, \eg, \cite[Sec.~6]{Andrews-Askey-Roy_2000})
that $\{\mathcal{L}_n^{\mu}\}_{n\in\Nat}$
is a complete orthogonal set in 
$
  \sii((0,\infty), x^\mu e^{-x} dx)
$. 
\end{proof}
\begin{Remark}
Note that the spectrum of~$L_\infty$ is doubly degenerate 
if $\Phi_B \in (\frac{1}{2}\Int)\setminus\Int$.
The situation $\Phi_B \in \Int$ corresponds to the classical 
partial wave decomposition, 
when all the eigenvalues of~$K$ are doubly degenerate except for one.
Finally, the spectrum of~$L_\infty$ is simple 
in the generic case $\Phi_B \not\in \frac{1}{2}\Int$.
\end{Remark}

\subsubsection{The strong-resolvent convergence}
Now we are in a position to prove the main auxiliary result of this paper. 
\begin{Proposition}\label{Prop.strong}
Let $\Phi_B \not\in \Int$.
Then the operator~$L_s$ converges to~$L_\infty$ 
in the strong-resolvent sense as $s \to \infty$, \ie,
\begin{equation*}
  \forall F \in \sii \,, \qquad
  \lim_{s \to \infty}
  \left\|
  L_s^{-1} F - L_\infty^{-1} F
  \right\|
  = 0
  \,.
\end{equation*}
\end{Proposition}
\begin{proof}
It follows from~\eqref{diamagnetic.bound} that~$0$ belongs 
to the resolvent set of~$L_s$ for all $s \geq 0$. 
For any fixed $F \in \sii$, let us set $\psi_s := L_s^{-1}F$.
In other words, $\psi_s$~satisfies the resolvent equation
\begin{equation}\label{re}
  \forall v \in \Dom(l_s) \,, \qquad
  l_s(v,\psi_s)
  = (v,F)
  \,,
\end{equation}
where~$l_s(\cdot,\cdot)$ is the sesquilinear form associated with~\eqref{J0.form}
and~$(\cdot,\cdot)$ denotes the inner product in~$\sii$.
In particular, choosing~$\psi_s$ for the test function~$v$ in~\eqref{re},
we have
\begin{equation}\label{resolvent.identity}
  l_s[\psi]
  = (\psi_s,F)
  \leq \|\psi_s\| \|F\|
  \,.
\end{equation}

Noticing that~\eqref{diamagnetic.bound} 
yields a Poincar\'e-type inequality 
$
  l_s[\psi_s] \geq \frac{1}{2} \|\psi\|^2
$
for any $\psi\in\Dom(l_s)$, we get from~\eqref{resolvent.identity}
the uniform bound
\begin{equation}\label{bound0}
  \|\psi_s\| \leq 2 \, \|F\|
  \,.
\end{equation}
At the same time, expressing the form~$l_s$ in polar coordinates 
(\cf~\eqref{form.polar}), the bounds \eqref{resolvent.identity}
and~\eqref{bound0} yield
\begin{equation}\label{bounds}
  \|\partial_r\psi_s\|^2 \leq 2 \, \|F\|^2
  \,, \quad
  \|r^{-1}(\partial_\theta\psi_s-i\alpha_s\psi_s)\|^2 \leq 2 \, \|F\|^2
  \,, \quad
  \|r\psi_s\|^2 \leq 32 \, \|F\|^2
  \,,
\end{equation}
where $\alpha_s(r,\theta) := \alpha(e^{s/2} r,\theta)$. 

Now, to get a similar estimate on the angular derivative of~$\psi_s$,
we employ the Hardy-type inequality~\eqref{Hardy} 
due to Laptev and Weidl as follows.  
Defining a new function $u_s\in\sii$ by
$\psi_s(y)=e^{s/2} u_s(e^{s/2} y)$
(\cf~the self-similarity transformation~\eqref{SST}),
making the change of variables $x=e^{s/2} y$
and recalling~\eqref{sigma},
we have
\begin{align}\label{unself}
  l_s[\psi_s]
  &\geq  \int_{\Real^2} \big|
  \nabla\psi_s(y)-iA_s(y) \, \psi_s(y)
  \big|^2 \, dy 
  \nonumber \\
  &= e^s \int_{\Real^2} \big|
  \nabla u_s(x)-iA(x) \, u_s(x)
  \big|^2 \, dx 
  \nonumber \\
  &\geq e^s \, c \int_{\Real^2} \frac{|u_s(x)|^2}{1+|x|^2} \, dx
  \,,
  \nonumber \\
  &= c \int_{\Real^2} \frac{|\psi_s(y)|^2}{e^{-s}+|y|^2} \, dy
  \,.
\end{align}
Here the first inequality follows just by neglecting 
the harmonic-oscillator potential term in~\eqref{J0.form}.
The second inequality is the Hardy inequality~\eqref{Hardy}
(the hypotheses of~\cite[Thm.~1]{Laptev-Weidl_1999} are clearly 
satisfied for $B \in C_0^0(\Real^2)$ assumed in the present paper)
and holds with a positive~$c$ if $\Phi_B \not\in \Int$.
Using this bound together with~\eqref{resolvent.identity} and~\eqref{bound0}, 
we therefore get 
\begin{equation}\label{ri3}
  \int_{\Real^2} \frac{|\psi_s(y)|^2}{e^{-s}+|y|^2} \, dy
  \leq \frac{2}{c} \, \|F\|^2
  \,.
\end{equation}

As a consequence of~\eqref{ri3}, we get 
$$
  \left\|r^{-1}i\alpha_s\psi_s\right\|^2
  \leq \frac{2}{c} \, C \, \|F\|^2
  \qquad\mbox{with}\qquad
  C := \sup_{\rho\in(0,\infty), \,\theta\in[0,2\pi)} 
  \frac{1+\rho^2}{\rho^2} \, \alpha(\rho,\theta)^2 
  \,.
$$
The finiteness of the constant~$C$ follows from the asymptotic behaviour
of~$\alpha(\rho,\theta)$ for small and large~$\rho$;
indeed, $\alpha(\rho,\theta)$ equals 
the bounded function $\alpha_\infty(\theta)$
for all sufficiently large~$\rho$ and the bound 
$|\alpha(\rho,\theta)| \leq \rho^2 \|B\|_\infty$ holds 
for all $\rho \geq 0$ and $\theta\in[0,2\pi)$.
It means that we can conclude from the central inequality in~\eqref{bounds}
the desired uniform boundedness of the angular derivative  
\begin{equation}\label{bound.angular}
  \|r^{-1}\partial_\theta\psi_s\| \leq \tilde{C} \, \|F\|
  \,,
\end{equation}
where the constant~$\tilde{C}$ depends on~$c$ and~$C$.

It follows from~\eqref{bound0}, \eqref{bounds} and~\eqref{bound.angular}  
that $\{\psi_s\}_{s \geq 0}$ is a bounded family in the space 
$
  \Dom_0 :=
  W^{1,2}(\Real^2) \cap \sii\big(\Real^2, |y|^2 \, dy\big)
$
equipped with the norm $\|\cdot\|_{\Dom_0}$ given by
$
  \|\psi\|_{\Dom_0}^2
  := \|\nabla\psi\|^2+\|y\psi\|^2+\|\psi\|^2 
$.
Therefore it is precompact in the weak topology of~$\Dom_0$.
Let~$\psi_\infty$ be a weak limit point,
\ie, for an increasing sequence of positive numbers $\{s_n\}_{n\in\Nat}$
such that $s_n \to \infty$ as $n \to \infty$,
$\{\psi_{s_n}\}_{n\in\Nat}$ converges weakly to~$\psi_\infty$ in~$\Dom_0$.
Actually, we may assume that the sequence converges strongly in $\sii$
because~$\Dom_0$ is compactly embedded in $\sii$.

Since $\{\psi_{s_n}\}_{n\in\Nat}$ converges strongly to~$\psi_\infty$ 
in $\sii$  and the function 
$e^{-s}+|y|^2$ converges locally uniformly for $y \in \Real^2\setminus\{0\}$
to $|y|^2$ as $s\to\infty$,
we have 
$$
  \forall \phi \in C_0^\infty(\Real^2\setminus\{0\}) 
  \,, \qquad
  \bigg(\phi,\frac{\psi_{s_n}}{\sqrt{e^{-s}+|\cdot|^2}}\bigg)
  \xrightarrow[s\to\infty]{}
  \bigg(\phi,\frac{\psi_\infty}{|\cdot|}\bigg)
  \,.
$$
Since $C_0^\infty(\Real^2\setminus\{0\})$ 
is dense in $\sii(\Real^2\setminus\{0\})=\sii(\Real^2)$
and the uniform bound~\eqref{ri3} holds true,
we may conclude that
$$
  \frac{\psi_{s_n}}{\sqrt{e^{-s}+|\cdot|^2}}
  \xrightarrow[s\to\infty]{w}
  \frac{\psi_\infty}{|\cdot|}
  \qquad\mbox{in}\qquad 
  \sii
  \,.
$$
In particular, this means that $\psi_\infty\in \sii(\Real^2,|y|^{-2}\,dy)$.

Summing up, we have proved that
$\{\psi_{s_n}\}_{n\in\Nat}$ converges weakly to~$\psi_\infty$ 
in the space
$$
  \Dom_\infty :=
  W^{1,2}(\Real^2) \cap \sii\big(\Real^2, |y|^2 \, dy\big)
  \cap \sii\big(\Real^2, |y|^{-2} \, dy\big)
$$
equipped with the norm $\|\cdot\|_{\Dom_\infty}$ given by
$
  \|\psi\|_{\Dom_\infty}^2
  := \|\psi\|_{\Dom_0}^2+\||y|^{-1}\psi\|^2
$.
In view of~\eqref{topology}, $\Dom_\infty$~coincides 
with~$\Dom(l_\infty)$ as a vector space,
so that we know that~$\psi_\infty$ belongs 
to the form domain of the operator~$L_\infty$. 
Taking $v \in C_0^\infty(\Real^2\setminus\{0\})$ 
as the test function in~\eqref{re}, with~$s$ being replaced by~$s_n$,
and sending~$n$ to infinity, we then easily check 
with help of~\eqref{AB-limit-potential} that
\begin{equation*}
  l_\infty(v,\psi_\infty)
  = (v,F)
  \,.
\end{equation*}
Since, by definition~\eqref{J0.form.AB}, 
$C_0^\infty(\Real^2\setminus\{0\})$ is a core for~$l_\infty$, 
we conclude that
$\psi_\infty = L_\infty^{-1} F$,
for \emph{any} weak limit point of $\{\psi_s\}_{s \geq 0}$.
The proof is concluded by recalling that we also have
the strong convergence of~$\{\psi_{s_n}\}_{n\in\Nat}$ to~$\psi_\infty$ 
in~$\sii$ as $n \to \infty$.
\end{proof}
\begin{Remark}\label{Rem.unself}
The crucial step in the proof is certainly the usage
of the Hardy inequality in the second inequality of~\eqref{unself}.
Indeed, it enables one to estimate the angular derivative 
from the central term in~\eqref{bounds}
by controlling the function
$r^{-1} i\alpha_s\psi_s$ via a uniform bound in~$s$.
\end{Remark}

\subsubsection{Spectral consequences}
In general, the strong-resolvent convergence of Proposition~\ref{Prop.strong}
is not sufficient to guarantee the convergence of spectra.
However, in our case, since the spectra are purely discrete
(\ie, the operators~$L_s$ and~$L_\infty$  have compact resolvents),
the eigenprojections converge even in norm (\cf~\cite{Weidmann_1980}).

In particular, $\lambda(s)$~converges to the lowest eigenvalue
of~$L_\infty$ as $s \to \infty$.
Recalling Proposition~\ref{Prop.AB-spectrum}
and the definition of~$\beta$ in~\eqref{Hardy.AB},
together with Proposition~\ref{Prop.sharp} 
for the case $\Phi_B\in\Int$, 
we thus conclude this long subsection by the ultimate result:
\begin{Corollary}\label{Corol.strong}
One has
$$
  \lambda(\infty) := \lim_{s\to\infty} \lambda(s) = (1+\beta)/2
  \,.
$$
\end{Corollary}
%

\subsection{A lower bound to the decay rate}\label{Sec.improved}
%
We come back to~\eqref{spectral.reduction.integral}.
It follows from Corollary~\ref{Corol.strong}
that for arbitrarily small positive number~$\eps$
there exists a (large) positive time~$s_\eps$ such that
for all $s \geq s_\eps$, we have $\lambda(s) \geq \lambda(\infty) - \eps$.
Hence, fixing $\eps>0$, for all $s \geq s_\eps$, we have
\begin{align*}
  {-\int_0^s \lambda(\tau) \, d\tau}
  &\leq {-\int_0^{s_\eps} \lambda(\tau) \, d\tau} {-[\lambda(\infty)-\eps](s-{s_\eps})}
  \\
  &\leq {[\lambda(\infty)-\eps] s_\eps} {-[\lambda(\infty)-\eps] s}
  \,,
\end{align*}
where the second inequality is due to the fact
that~$\lambda(s)$ is non-negative for all $s \geq 0$
(it is in fact greater than or equal to~$1/2$,
\cf~\eqref{diamagnetic.bound}).
At the same time, assuming $\eps \leq 1/2$, we trivially have
$$
  {-\int_0^s \lambda(\tau) \, d\tau}
  \leq 0
  \leq {[\lambda(\infty)-\eps] s_\eps} {-[\lambda(\infty)-\eps] s}
$$
also for all $s \leq s_\eps$.
Summing up, for every $s \in [0,\infty)$, we have
\begin{equation}\label{instead}
  \|\tilde{u}(s)\|_{K}
  \leq C_\eps \, e^{-[\lambda(\infty)-\eps]s} \, \|\tilde{u}_0\|_{K}
  \,,
\end{equation}
where $C_\eps := e^{s_\eps} \geq e^{[\lambda(\infty)-\eps]s_\eps}$.

Now we return to the original variables $(x,t)$ via~\eqref{space-times}.
Using~\eqref{preserve} together with the point-wise estimate $1 \leq K$,
and recalling that $\tilde{u}_0=u_0$,
it follows from~\eqref{instead} that
$$
  \|u(t)\|
  = \|\tilde{u}(s)\|
  \leq \|\tilde{u}(s)\|_{K}
  \leq C_\eps \, (1+t)^{-[\lambda(\infty)-\eps]} \, \|u_0\|_{K}
$$
for every $t \in [0,\infty)$.
Consequently, we conclude with
$$
  \|e^{-H_B t}\|_{\sii(K) \to \sii}
  \equiv \sup_{u_0 \in \sii(K)\setminus\{0\}}
  \frac{\|u(t)\|}{\ \|u_0\|_{K}}
  \leq C_\eps \, (1+t)^{-[\lambda(\infty)-\eps]}
$$
for every $t \in [0,\infty)$.
Since~$\eps$ can be made arbitrarily small,
this bound implies
$$
  \gamma_B \geq \lambda(\infty) = (1+\beta)/2
  \,.
$$
This together with Corollary~\ref{Corol.straight}
proves Theorem~\ref{Thm.rate}.

\subsection{A global upper bound to the heat semigroup}\label{Sec.alt}
%
Theorem~\ref{Thm.rate} provides quite precise information
about the extra polynomial decay of solutions~$u$ of~\eqref{heat.intro}
in the sense that the decay rate~$\gamma_B$
is  better by a factor~$\beta/2$ with respect to the field-free case.
On the other hand, we have no control over
the constant~$C_\gamma$ in~\eqref{rate}
(in principle it may blow up as $\gamma \to \gamma_B$).
We therefore conclude this section by establishing
a global (in time) upper bound to the heat semigroup
(\ie~we get rid of the constant~$C_\gamma$)
but the prize we pay is just a qualitative knowledge
about the decay rate.
It is a consequence of~\eqref{spectral.reduction.integral}
and the following improvement upon~\eqref{diamagnetic.bound}:
\begin{Proposition}\label{Prop.positivity}
Let $B\not=0$. Then
$$
  \forall s \in [0, \infty)
  \,,\qquad
  \lambda(s) > 1/2
  \,.
$$
\end{Proposition}
\begin{proof}
The fact that $\lambda(s) \geq 1/2$ for all $s \geq 0$,
regardless of the presence of~$B$,
is due to~\eqref{diamagnetic.bound}.
To show that the inequality is strict,
let us assume by contradiction that $B\not=0$
and that $\lambda(s) = 1/2$ for some $s \geq 0$.
Let~$v$ denote a corresponding eigenfunction of~$L_s$.
We have
\begin{equation}\label{identities}
  0
  = l_s[v] - \frac{1}{2} \, \|v\|^2
  \geq \big\|\nabla|v|\big\|^2 + \frac{1}{16} \, \|y v\|^2 - \frac{1}{2} \, \|v\|^2
  \geq 0
  \,,
\end{equation}
where the first estimate is due the diamagnetic inequality~\eqref{diamagnetic}
and the second inequality follows from the variational characterization
of the lowest eigenvalue of~$L_\mathrm{HO}$ (which is $1/2$).
Hence, $|v|$~must coincide with the ground state of~$L_\mathrm{HO}$.
The latter is unique (up to a normalization factor)
and can be chosen positive.
Consequently, $v$~as an eigenfunction of~$L_s$ is unique
(up to a normalization factor) and can be chosen positive.
The equality in~\eqref{identities} then implies
$$
  \|A_s v\| = 0
  \,.
$$
Hence, employing the continuity of~$A_s$ and~$v$, $A_s=0$ identically,
which in turn implies that $B$~is trivial, too.
\end{proof}
\begin{Remark}
If $B\not=0$ but $\Phi_B\in\Int$, then it follows from
Proposition~\ref{Prop.sharp} that $\lambda(s)$ converges to~$1/2$ as $s\to\infty$,
despite of the strict inequality of Proposition~\ref{Prop.positivity}.
\end{Remark}

Combining Proposition~\ref{Prop.positivity} with Corollary~\ref{Corol.strong},
we see that the number
\begin{equation}
  c_B := \inf_{s \in [0,\infty)}\lambda(s) - 1/2
\end{equation}
is positive if, and only if, $\Phi_B \not\in \Int$.
In any case, \eqref{spectral.reduction.integral}~implies
$$
  \|\tilde{u}(s)\|_{K}
  \leq \|\tilde{u}_0\|_{K} \, e^{-(c_B+1/2)s}
$$
for every $s \in [0,\infty)$.
Using this estimate instead of~\eqref{instead},
but following the same type of arguments as in Section~\ref{Sec.improved}
below~\eqref{instead}, we thus conclude with:
\begin{Theorem}\label{Thm.rate.alt}
Let $B \in C_0^0(\Real^2)$.
We have
$$
  \forall t \in [0,\infty), \qquad
  \big\|e^{-H_B t}\big\|_{\sii(K) \to \sii}
  \,\leq\, (1+t)^{-(c_B+1/2)}
  \,,
$$
where $c_B > 0$ if $\Phi_B \not\in \Int$
(and $c_B = 0$ if $\Phi_B \in \Int$).
\end{Theorem}
%

\section{Conclusions}\label{Sec.end}
%
The results of this paper show that the presence
of a local magnetic field leads to an improvement of the large-time
decay of solutions of the heat equation.
As explained in the introduction,
the physical reason behind the better decay
is the diamagnetic effect of the magnetic field.
This is mathematically expressed in terms of the existence
of Hardy-type inequalities for the magnetic Schr\"odinger operator,
which played a central role in our proofs.

It follows from our analysis that the solutions
of the heat equation for large time behave as
if the magnetic field were replaced by
the singular Aharonov-Bohm magnetic field with the same total flux.
This also explains the role of the ``degeneracies'' $\Phi_B\in\Int$,
since they correspond to the flux quanta for which the Aharonov-Bohm effect 
is not detectable in our setting (\cf~\cite[Sec.~10.2.4]{Griffiths}). 

Since no assumptions about the symmetry of the magnetic field
are made whatsoever in the present paper, it can be also viewed
as a generalization of some of the results of the recent work~\cite{Kovarik_2010}

It is known that the method of self-similarity variables
yield sharp decay rates.
Therefore, we conjecture that the inequality of Theorem~\ref{Thm.rate}
can be replaced by equality, \ie, $\gamma_B=(1+\beta)/2$ for any field~$B$.
A careful analysis in the spirit of \cite{Escobedo-Zuazua_1991,Duro-Zuazua_1999}
could even show that the solutions can indeed be approximated
by the solutions of the effective equation with the Aharonov-Bohm field.

The question of optimal value of the constant~$c_B$
(and its quantitative dependence on the magnetic field~$B$,
especially on the total magnetic flux~$\Phi_B$)
from Theorem~\ref{Thm.rate.alt} also constitutes an interesting open problem.
Note that the constant is related to the decay rate 
by $c_B+1/2 \leq \gamma_B$.

We expect the same decay rates if the assumption about the compact
support of~$B$ is replaced by its fast decay at infinity.
However, it is quite possible that a slow decay of the field
at infinity will improve the decay of the solutions even further.
In particular, can~$\gamma_B$ be strictly greater than~$(1+\beta)/2$
if~$B$ decays to zero very slowly at infinity?

The characteristic hypothesis of the present paper, under which
the additional decay rate is proved, is that the total magnetic
flux is not an integer. Even if this hypothesis is violated
and the magnetic field is non-trivial, however, there exist
Hardy-type inequalities for the magnetic Hamiltonian~\cite{Weidl_1999}.
Moreover, it is shown in~\cite{Kovarik_2010} that there is
a logarithmic improvement to the decay of the heat kernel
even if the total flux is an integer (and the field is rotationally symmetric).
It should be possible to establish the additional
logarithmic decay rate by the method of the present paper
(without any symmetry assumption),
by simply computing the next term in the asymptotic expansion
of the eigenvalue $\lambda(s)$ as $s\to\infty$.

More generally, recall that we expect that there is always
an improvement of the decay rate for the heat semigroup
of an operator satisfying a Hardy-type inequality
(\cf~\cite[Conjecture in Sec.~6]{KZ1} and~\cite[Conjecture~1]{FKP}).
The present paper confirms the general conjecture in the particular case
of two-dimensional magnetic Schr\"odinger operators.

Finally, it would be interesting to examine the effect 
of the presence of the local two-dimensional magnetic field 
in other physical models. 
As one possible direction of this research, let us mention 
the problem of Strichartz estimates for
the magnetic Schr\"odinger equation studied recently in higher dimensions
in \cite{D'Ancona-Fanelli-Vega-Visciglia_2010}.

\section*{Acknowledgement}
The author would like to thank Carlos Mora-Corral for useful discussions.
The work was partially supported by the Czech Ministry of Education,
Youth and Sports within the project LC06002
and by the GACR grant P203/11/0701.

%
{\small
\providecommand{\bysame}{\leavevmode\hbox to3em{\hrulefill}\thinspace}
\providecommand{\MR}{\relax\ifhmode\unskip\space\fi MR }
\providecommand{\MRhref}[2]{%
  \href{http://www.ams.org/mathscinet-getitem?mr=#1}{#2}
}
\providecommand{\href}[2]{#2}

%
}
\end{document}